\documentclass{article}

\usepackage{amsfonts}
\usepackage{amsthm}

\newtheorem{thm}{Theorem}[section]
\newtheorem{prop}[thm]{Proposition}
\newtheorem{lem}[thm]{Lemma}

\newtheorem{remark}[thm]{Remark}
\newtheorem{defn}[thm]{Definition}

\title{Coincident root loci and Jack and Macdonald polynomials for special
values of the parameters}
\author{M. Kasatani \thanks {Department of Mathematics, Graduate School of Science, 
Kyoto University, Kyoto 606-8502, Japan, e-mail: {\tt kasatani@math.kyoto-u.ac.jp}}, T. Miwa \thanks {Department of Mathematics, Graduate School of Science, Kyoto University, Kyoto 606-8502, Japan, e-mail: {\tt tetsuji@math.kyoto-u.ac.jp}}, A. N. Sergeev \thanks {Balakovo Institute of Technology and
Control, Balakovo, 413800, Russia, e-mail: {\tt sergeev@bittu.org.ru}}, A.P. Veselov \thanks {Department of Mathematical Sciences, Loughborough University, Loughborough,  LE11 3TU, UK and  
Landau Institute for Theoretical Physics, Kosygina 2, Moscow, 117940, Russia, e-mail: {\tt A.P.Veselov@lboro.ac.uk} }}
\date{\empty}

\begin{document}
\maketitle

{\small  {\bf Abstract.} We consider the coincident root loci
consisting of the polynomials with at least two double roots and
present a linear basis of the corresponding ideal in the algebra of symmetric polynomials 
in terms of the
Jack polynomials with special value of parameter $\alpha = -2.$ As
a corollary we present an explicit formula for the Hilbert-Poincar\`e series of this ideal and 
the generator of the minimal degree as a special Jack polynomial.

A generalization to the case of the symmetric polynomials
vanishing on the double shifted diagonals
and the Macdonald polynomials specialized at $t^2 q = 1$ is also
presented. We also give similar results for the interpolation
Jack polynomials.}

\section{Introduction}
In 1857 Arthur Cayley published a short paper \cite{Cayley} where
he considered a problem (which he prescribed to Sylvester) of how
to determine when a polynomial of degree $n$ has a multiple root
of multiplicity at least $m.$ When $m=2$ the answer is of course
well known: the corresponding algebraic variety is called the {\it
discriminant} and can be defined by equating the discriminant
of a polynomial to zero. The case of general $m$ corresponds to
the natural strata in the discriminant also known as {\it
coincident root loci}.

Cayley considered also a more general question when a polynomial
has several multiple roots with prescribed multiplicities. One can label the corresponding stratum in the discriminant  by a partition $\mu = (\mu_1, \ldots, \mu_k), \mu_1 + \ldots + \mu_k = n$. 
Cayley suggested an approach to this problem based on the classical invariant theory of the binary forms  and demonstrated it for quartics and quintics. 

For general $n$ and one multiple root (i.e. for the partition $\mu = (m, 1, 1, \dots, 1) = (m, 1^{n-m})$)
some interesting results were found by J. Weyman
\cite{W1,W2}
(see also recent paper \cite{Chip}), but the problem is still largely open even in this case.
In terms of the symmetric functions of the roots it can be formulated as
follows:
describe the ideal $I_n(m)$ in the algebra of symmetric polynomials of $n$
variables
vanishing when $m$ of the variables are equal. The ring of functions on the
the corresponding coincident root loci is isomorphic to ${\bf
C}[x_1,\ldots,x_n]^{\mathfrak{S}_n}/I_n(m)$.

Recently Feigin, Jimbo, Miwa and Mukhin \cite{FJMM} described
explicitly a linear basis in this ideal in terms of the Jack
polynomials with special value of the parameter and special
combinatorics of the corresponding Young diagrams. Some
explanation for the appearance of the Jack polynomials in this
problem was found by Sergeev and Veselov \cite{SV}, who came to a
similar problem in a different way investigating the deformed
Calogero-Moser operators. In particular, the results of  \cite{SV}
suggest that one should be able to generate the ideals by certain
Jack polynomials also for the partitions
$\mu = (p^k, q^{l})$, i.e. in the case when we have $k$ roots of
multiplicity $p$ and $l$ roots of multiplicity $q$, $p k + q l = n$.

In the present paper we show that this is true in the simplest case
when we have two double roots: $\mu = (2^2,1^{n-4})$.
More precisely, consider the ideal $I_n$ consisting of
 the symmetric polynomials of $n$ variables $P(x_1,x_2, \dots, x_n)$ with
the property $P(x_1,x_2, \dots, x_n)=0$ if $x_1 = x_2$ and $x_3 = x_4.$
We construct explicitly a linear basis of this ideal in terms of the Jack
polynomials with the special value of the parameter $\alpha=-2$
(in Macdonald's notation \cite{Mac}). An interesting novelty in comparison
with \cite{FJMM} is
that for some diagrams the Jack polynomial itself is not good enough and one
should
consider a certain linear combination of two Jack polynomials with different
Young diagrams
(see Section 2,3 for details).

We consider also two related ideals $J_n$ and $I_n^*$ consisting of
symmetric polynomials vanishing when  $x_1 = t x_2, x_3 =  t x_4$ and $x_2 =
x_1 + 1/2, x_4 = x_3 + 1/2$ respectively
and construct explicit linear bases in terms of the Macdonald polynomials
specialized at $q = t^{-2}$
(for $J_n$) and in terms of the interpolation Jack polynomials introduced
and investigated in \cite{KS,OO}.

The structure of the paper is following. First we introduce the
admissible partitions and construct a basis in the ideal $J_n$ in
terms of the corresponding Macdonald polynomials. The proof is
based on the results and ideas from \cite{FJMM,FS,FJMM02,K}.
Then we show how to derive from this the description of the ideals
$I_n$ and $I^*_n$ in terms of Jack polynomials and their
generalizations introduced by Knop, Sahi, Okounkov and Olshanski
\cite{KS,OO}. 

In the last section as a corollary of our results we prove the following explicit formula 
for the Hilbert-Poincar\`e series of the ideal $I_n$:
\begin{eqnarray}
{\mathop{\mathrm{ch}}\nolimits}\,
I_n:=\sum_{d=0}^{\infty} (\dim I_{n,d}) q^d = \frac{q^{(n-1)(n-3)}}{(q)_1(q)_{n-3}}
+\frac{q^{(n-2)(n-1)}}{(q)_1(q)_{n-2}}
+\frac{q^{n(n-1)}}{(q)_n} \nonumber
\end{eqnarray}
where $(q)_s = \prod_{j=1}^{s}(1-q^j).$ 
We present also a generator of the minimal degree as a special Jack polynomial.

\section{Symmetric polynomials and zero condition on the double diagonals}
Throughout this paper we assume that $n\geq 4$. Let
$\Lambda_n=\mathbb{C}[x_1,\cdots,x_n]^{{\mathfrak S}_n}$ be the
ring of symmetric polynomials. We extend this space to
$\Lambda_n\otimes_\mathbb{C}\mathbb{C}(t)$, and consider the zero
condition on the double "$t$-diagonals":
\begin{equation}
f(x,tx,y,ty,x_5,\cdots,x_n)=0.\label{WC}
\end{equation}
We denote by $J_n$ the subspace of
$\Lambda_n\otimes_\mathbb{C}\mathbb{C}(t)$
consisting of the polynomials satisfying (\ref{WC}), and
by $J_{n,d}$ the degree $d$ component of $J_n$.

We will construct a basis of $J_n$ by using the Macdonald
polynomials $P_\lambda$ specialized at $q=t^{-2}$, where $\lambda$
belongs to a certain set of partitions of $n$. In \cite{FJMM02}, a
similar result is established for the zero condition on a single
$t$-diagonal of higher codimensions:
\begin{equation}\label{K}
z_i/z_{i+1}=t\quad(1\leq i\leq k)
\end{equation}
by using the specialization of the Macdonald polynomials at $q=t^{-k-1}$.
In fact, in \cite{FJMM02}, a more general result for the case of
several distinct shifted-diagonals is obtained.
The condition (\ref{WC}) is related to the case $k=1$ in (\ref{K}),
which is trivial because the condition is equivalent to that
$f$ is divisible by the square of the discriminant. In this paper,
we study the new case of double shifted-diagonals.

Let us introduce the set of admissible partitions which is relevant
to the present case. We use the dual operator language (see \cite{FS}).
Denote by $R$ the polynomial ring $\mathbb{C}[e_0,e_1,e_2,\ldots]$.
We count weight of $e_i$ as $1$ and degree of $e_i$ as $i$, and
denote by $R_{n,d} \subset R$ the subset consisting of the weight $n$ and
degree $d$ polynomials.
\begin{defn}\normalfont
A monomial
\begin{equation}\label{MONO}
m=e_0^{a_0}e_1^{a_1}e_2^{a_2}\cdots
\end{equation}
is called non-admissible if and only if $m=m_0m'$ ($m_0,m'\in R$) and
one of the following is valid:
\begin{eqnarray}
&&\hbox{$m_0=m_1m_2$,  where $m_1$ and $m_2$ are of the form $e_i^2$
or $e_je_{j+1}$}, \quad \label{N1}\\
&&m_0=e_i^3e_{i+2}.\label{N2}
\end{eqnarray}
In (\ref{N1}) $i$ is not necessarily distinct from $j,j+1$.
A monomial is called admissible if and only if it is not non-admissible.
\end{defn}

We denote by $\pi_n$ the set of partitions of length $n$,
$\lambda=(\lambda_1,\lambda_2,\cdots,\lambda_n)$.
We define the degree $d$ component
$\pi_{n,d}=\{\lambda\in \pi_n; |\lambda|=d\}$,
where $|\lambda|=\lambda_1+\cdots+\lambda_n$.
We set $e_\lambda=\prod_{i=1}^n e_{\lambda_i}$.

There is a one-to-one correspondence between a monomial $m\in R_{n,d}$
given by (\ref{MONO}) and a partition $\lambda\in\pi_{n,d}$:
\begin{equation}
a_i=\sharp\{j;\lambda_j=i\}.
\end{equation}
We say a partition $\lambda$ is admissible if and only if the corresponding
monomial $m$ is admissible.
Admissible partitions are classified into three different cases:

\medskip
{\it Case A}\quad There exists $i$ such that
\begin{eqnarray}
&&\lambda_i=\lambda_{i+1}=\lambda_{i+2}=\lambda_{i+3}+2,\\
&&\lambda_{i-1}-\lambda_i > 2 \mbox{\quad if $i \geq 2$},\\
&&\lambda_j-\lambda_{j+1}\geq 2\mbox{\quad if $j\leq i-2$ or $j\geq i+3$}.
\end{eqnarray}

{\it Case B}\quad There exists $i$ such that
\begin{eqnarray}
&&\lambda_i=\lambda_{i+1}+1=\lambda_{i+2}+2,\\
&&\lambda_j-\lambda_{j+1}\geq 2\mbox{\quad if $j\leq i-1$ or $j\geq i+2$}.
\end{eqnarray}

{\it Case C}\quad
$\lambda$ is admissible and belongs to neither {\it Case A} nor {\it B}.

\medskip
Now, we consider the ring of symmetric polynomials
$\Lambda_n\otimes_\mathbb{C}\mathbb{C}(q,t)$ with coefficients in
$\mathbb{C}(q,t)$.
Following \cite{Mac}, we define a homomorphism

\[
u_\lambda:\Lambda_n\otimes_\mathbb{C}\mathbb{C}(q,t)
\rightarrow\mathbb{C}(q,t)
\]
by
\begin{equation}
u_\lambda(f(x_1,\cdots,x_n))=
f(t^{n-1}q^{\lambda_1},t^{n-2}q^{\lambda_2},\cdots,q^{\lambda_n}).
\end{equation}

The Macdonald polynomial $P_\lambda$ corresponding to a partition $\lambda$
is an eigenvector of a difference operator $D_n(X;q,t)$ (see \cite{Mac}).
The corresponding eigenvalue is given by
\begin{equation}
\prod_{i=1}^n (1+Xt^{n-i}q^{\lambda_i}).
\end{equation}
We define an element $E(\lambda)\in\mathbb{C}(q,t)^n/\mathfrak{S}_n$ by
\begin{equation}
E(\lambda)=\{t^{n-i}q^{\lambda_i}\}_{1\leq i\leq n}.
\end{equation}

When we specialize at $q=t^{-2}$, it may happen that
$E(\lambda)=E(\nu)$ for $\lambda\not=\nu$. In fact, in {\it Case A} or {\it
B},
we show in the below that
there exists $\nu\not=\lambda$ such that $E(\lambda)=E(\nu)$ for $q=t^{-2}$.
We will show that if a partition $\lambda$ belongs to {\it Case C},
the Macdonald polynomial $P_\lambda$ is well-defined at $q=t^{-2}$
and satisfies the condition (\ref{WC}). A similar statement is not true in
{\it Case A} or {\it B}. We define a modified polynomial
$\bar{P}_\lambda$ to overcome this difficulty.

In {\it Case A}, we set
\begin{eqnarray}\label{BAR}
\bar{P}_\lambda=u_0(P_\lambda)
\left(\frac{P_\lambda}{u_0(P_\lambda)}-\frac{P_\nu}{u_0(P_\nu)}\right)
\end{eqnarray}
where the partition $\nu$ is given by
\begin{eqnarray}
&&(\nu_i,\nu_{i+1},\nu_{i+2},\nu_{i+3})
=(\lambda_i+1,\lambda_{i+1}-1,\lambda_{i+2}-1,\lambda_{i+3}+1),\\
&&\nu_j=\lambda_j \quad(j\leq i-1 \mbox{ or } j\geq i+4).
\end{eqnarray}
In {\it Case B}, we define $\bar{P}_\lambda$ by (\ref{BAR}) where $\nu$
is given by
\begin{eqnarray}
&&(\nu_i,\nu_{i+1},\nu_{i+2})=(\lambda_i-1,\lambda_{i+1},\lambda_{i+2}+1),\\
&&\nu_j=\lambda_j \quad(j\leq i-1 \mbox{ or } j\geq i+3).
\end{eqnarray}
For convenience, we set $\bar{P}_\lambda=P_\lambda$ in {\it Case C}.
Note that in {\it Case A} or {\it B}, we have $E(\lambda)=E(\nu)$.

The main result of this paper is
\begin{thm}\label{thm:main}
The modified polynomial $\bar{P}_\lambda$ has no pole at $t^2q=1$ and
the specialization $\bar{P}_\lambda|_{q=t^{-2}}$
satisfies the zero condition $(\ref{WC})$.
The set of polynomials $\{\bar{P}_\lambda|_{q=t^{-2}};\lambda\in\pi_{n}
\mbox{ is admissible}\}$ is a basis of $J_n$.
\end{thm}

In the limit $t\rightarrow1$ the condition (\ref{WC}) becomes
\begin{equation}
f(x,x,y,y,x_5,\cdots,x_n)=0 \label{xxyy}
\end{equation}
which defines the ideal $I_n.$ As a corollary we have a similar
claim for this ideal and the corresponding modified Jack
polynomials specialized at $\alpha = -2$ (see Section 5 for the
details).

\begin{remark}\label{rem:construction}\normalfont
For the partition $\lambda$ of Case A,
 there does not exist a partition $\mu<\lambda$ such that $E(\lambda)=E(\mu)$.
The same statement holds for $\nu$.
(The proof is given in Lemma \ref{lem:caseA}.
Note that $E(\lambda)=E(\nu)$. However,
 neither $\nu<\lambda$ nor $\lambda<\nu$ is valid.)
Hence, from Lemma \ref{lem:no pole}, we have
 the Macdonald polynomials $P_\lambda$ and $P_\nu$ have no pole at $t^2q=1$,
 and they belong to the same eigenspace for
 the Macdonald operator $D_n(X;t^{-2},t)$.
However, these polynomials do not safisfy the condition (\ref{WC}).
This is why we consider $\bar{P}_\lambda$ in this case.

In Case B, the Macdonald polynomial $P_\nu$ has no pole.
 This is because $\nu$ belongs to Case C.
The coefficient $\frac{u_0(P_\lambda)}{u_0(P_\nu)}$ has a single pole
 at $q=t^{-2}$. (See the formula (\ref{SPEC}) in Section \ref{CONST}.)
Therefore, the Macdonald polynomial $P_\lambda$ has a single pole at $q=t^{-2}$.
Theorem \ref{thm:main} asserts that the modified polynomial $\bar{P}_\lambda$
 has no pole, and moreover it satisfies the condition (\ref{WC}).

The above argument implies that the space $J_n$ is invariant under
 the action of the Macdonald operator (cf. \cite{SV2}, where a more general fact is proved).
In Case A and C, $\bar{P}_\lambda$ is an eigenfunction of the operator.
In Case B, $\bar{P}_\lambda$ is not an eigenfunction, but
 $\bar{P}_\lambda$ and the eigenfunction $\bar{P}_\nu=P_\nu$
 constitute a Jordan block.
\end{remark}

\begin{lem}\label{lem:caseA}
Let $\lambda$ be any partition of Case A.
Fix the integer $i$ such that
 $\lambda_{i}=\lambda_{i+1}=\lambda_{i+2}=\lambda_{i+3}+2$.
Then there exist only two patitions $\mu$ $(\neq\lambda)$
 such that $E(\lambda)=E(\mu)$.
Precisely, $(\mu_{i},\mu_{i+1},\mu_{i+2},\mu_{i+3})=
(\lambda_{i}+1,\lambda_{i}-1,\lambda_{i}-1,\lambda_{i}-1)$
 or $(\lambda_{i}+1,\lambda_{i},\lambda_{i}-1,\lambda_{i}-2)$.
(All other components of $\mu$ are equal to those of $\lambda$.)
\end{lem}
\begin{proof}
Suppose that there exists $\mu$ such that $E(\lambda)=E(\mu)$ and $\lambda\neq\mu$.
When we specialize $q=t^{-2}$, we can identify $E(\lambda)$ with
 the set of integers $I(\lambda):=\{j+2\lambda_j; j=1,\cdots,n\}$.
We separate $I(\lambda)$ to odd and even parts:
Define $I_1(\lambda):=\{j+2\lambda_j; j=i,i\pm 2,i\pm 4,\cdots \}$ and $I_2(\lambda):=I(\lambda) \backslash I_1(\lambda)$.
Then, $I_1(\lambda)=I_1(\mu)$ and $I_2(\lambda)=I_2(\mu)$.

We consider the part $I_1(\lambda)$.
Note that for $p,s\in i+2\mathbb{Z}$,
\begin{eqnarray}
&&\mbox{$\lambda_p-\lambda_{s}\geq s-p$ \quad if $s-p\geq 4$}, \label{data1-1}\\
&&\mbox{$\lambda_p-\lambda_{s}> s-p$ \quad if $p \neq i$ and $s-p\geq 2$}, \label{data1-2}\\
&&\mbox{$\lambda_i-\lambda_{s}> s-i$ \quad if $s-i\geq 6$}, \label{data1-3}
\end{eqnarray}

Suppose that there exists $j_1\in i+2\mathbb{Z}$ such that $\lambda_{j_1} \neq \mu_{j_1}$ and fix the minimum one.
Then there exists $j_2>j_1$ such that $j_1+2\lambda_{j_1}=j_2+2\mu_{j_2}$,
and there exists $j_3 \geq j_1$ such that $j_3 \neq j_2$ and $j_2+2\lambda_{j_2}=j_3+2\mu_{j_3}$.
Inductively, we define $j_\alpha \geq j_1$ such that $j_\alpha \neq j_{\alpha-1}$ and $j_{\alpha-1}+2\lambda_{j_{\alpha-1}}=j_\alpha+2\mu_{j_\alpha}$ for $\alpha\geq 4$.
For $j_2$ and $j_3$, there are two cases: $j_2<j_3$ or $j_2>j_3$.
We will show by induction that $j_2<j_3 \Rightarrow j_1<j_2<j_3<j_4< \cdots$.
Hence the condition $j_2<j_3$ contradicts the finiteness of indexes.

Let $\alpha\geq 4$ and suppose that $j_1<j_2< \cdots < j_{\alpha-1}$.
To show $j_{\alpha-1}<j_\alpha$, we prove step by step $j_2<j_\alpha,j_3<j_\alpha,\cdots, j_{\alpha-1}<j_\alpha$.
Fix $\beta\leq \alpha-2$. If $\beta\geq 3$, suppose further that $j_{\beta-1}<j_\alpha$.
Then $j_{\alpha-1}-j_{\beta-1}\geq 4$.
Thus from (\ref{data1-1}),
\begin{eqnarray}
\mu_{j_\beta}-\mu_{j_\alpha}&=&\lambda_{j_{\beta-1}}-\lambda_{j_{\alpha-1}}
                      +(j_{\beta-1}+j_{\alpha}-j_{\beta}-j_{\alpha-1})/2 \nonumber\\
&>& (j_{\alpha-1}-j_{\beta-1})+j_{\beta-1}-j_{\alpha-1} \nonumber\\
&=& 0. \nonumber
\end{eqnarray}
Hence $j_\beta<j_\alpha$.
Inductively with respect to $\beta$, we have $j_{\alpha-2}<j_\alpha$.
If $j_{\alpha-1}-j_{\alpha-2}=2$, then since $j_{\alpha-2}<j_\alpha$ and $j_{\alpha-1} \neq j_\alpha$, we see $j_{\alpha-1}<j_\alpha$.
If $j_{\alpha-1}-j_{\alpha-2}\geq 4$, then from (\ref{data1-1}),
 the inequality $\mu_{j_{\alpha-1}}-\mu_{j_\alpha}>0$ holds. Hence $j_{\alpha-1}<j_\alpha$.

Next, we show that $j_1=i$ and $j_2-j_1=2$.

Assume that $j_1 \neq i$, then from (\ref{data1-2}),
 the inequality $\mu_{j_2}-\mu_{j_3}>0$ holds.
Hence $j_2<j_3$.
This leads to the contradiction.

Assume that $j_1=i$ and $j_2-j_1\geq 6$,
 then from (\ref{data1-3}), the inequality $\mu_{j_2}-\mu_{j_3}>0$ holds.
Hence $j_2<j_3$.
This leads to the contradiction.

Assume that $j_1=i$ and $j_2-j_1=4$.
If $j_3>j_1$, then from (\ref{data1-1}),
 we see the inequality $\mu_{j_2}-\mu_{j_3}>0$ holds.
This leads to the contradiction.
If $j_3=j_1$, repeat the above argument
 for the index set $\{i+2\mathbb{Z}\}\backslash\{i,i+4\}$.
Since this set does not contain the element $i$,
 we have $\lambda_{j}=\mu_{j}$ if $j\in \{i+2\mathbb{Z}\}\backslash\{i,i+4\}$.
However, $\mu_{i}-\mu_{i+2}=(\lambda_{i+4}+2)-\lambda_{i+2}<0$.
Hence it leads to the contradiction.

We have shown that $j_1=i$ and $j_2-j_1=2$.
Since $j_2>j_3$, we see $j_3=j_1$.
As stated above, $\lambda_{j} = \mu_{j}$ holds for $j\in \{i+2\mathbb{Z}\}\backslash\{i,i+2\}$.

Therefore, we have $(\mu_{i},\mu_{i+2})=(\lambda_{i}+1,\lambda_{i}-1)$ or $(\lambda_{i},\lambda_{i})$.

For the part $I_2(\lambda)$,
note that for $p,s\in \{i+1+2\mathbb{Z}\}$,
\begin{eqnarray}
&&\mbox{$\lambda_p-\lambda_{s}>s-p$ \quad if $s-p \geq 4$}, \label{data2-1}\\
&&\mbox{$\lambda_p-\lambda_{s}>s-p$ \quad if $p \neq i+1$ and $s-p \geq 2$}. \label{data2-2}
\end{eqnarray}
One can show that $(\mu_{i+1},\mu_{i+3})=(\lambda_{i+1}-1,\lambda_{i+1}-1)$ or $(\lambda_{i+1},\lambda_{i+1}-2)$.

Therefore, the only possiblility of the partition $\mu(\neq\lambda)$ is
\begin{eqnarray}
(\mu_{i},\mu_{i+1},\mu_{i+2},\mu_{i+3})&=&
(\lambda_{i}+1,\lambda_{i}-1,\lambda_{i}-1,\lambda_{i}-1) \nonumber\\
&\mbox{or}&(\lambda_{i}+1,\lambda_{i},\lambda_{i}-1,\lambda_{i}-2). \nonumber
\end{eqnarray}
\end{proof}

\section{Dimension of the space of polynomials $J_{n,d}$}
In this section we give an upper estimate for
$\dim_{\mathbb{C}(t)} J_{n,d}$. In Section \ref{CONST}, we will
show the estimate is exact.

We denote by $I_n$ the subspace of $\Lambda_n$ consisting of
the polynomials satisfying (\ref{xxyy}), and by $I_{n,d}$
its degree $d$ component.
Note that $\dim_{\mathbb{C}(t)}J_{n,d}\leq\dim_{\mathbb{C}}I_{n,d}$.

In the dual language, the above condition (\ref{xxyy}) is equivalent
to a quartic relation for the abelian current $e(x)=\sum_{i\geq0}e_ix^i$.
The relation is
\begin{equation}
e(x)^2e(y)^2=0.
\end{equation}
Let $r_{i,j}$ be the coefficients in $e(x)^2e(y)^2$:
\begin{equation}
e(x)^2e(y)^2=\sum_{i,j}r_{i,j}x^iy^j.
\end{equation}
We denote by $J^*$ the ideal of $R$ generated by $\{r_{i,j}\}$, and set
$J^*_n=J^*\cap R_n$. For $r\in R$ we denote by $\bar r$
the image of $r$ in $ R/J^*$.
\begin{prop}
The image of the admissible monomials $e_\lambda$ of weight $n$
spans the quotient space $R_n/J^*_n$.
\end{prop}

\begin{proof}
We introduce an ordering of monomials.
For $e_\lambda=e_0^{a_0}e_1^{a_1}e_2^{a_2}\cdots$ and
$e_\mu=e_0^{b_0}e_1^{b_1}e_2^{b_2}\cdots$, we write $e_\lambda \succ e_\mu$
if
and only if $e_\lambda$ and $e_\mu$ have the same weight and degree,
and $a_0<b_0$ or $a_0=b_0,a_1<b_1$ or $a_0=b_0,a_1=b_1,a_2<b_2$ or $\cdots$.
We also write $\lambda\succ\mu$. This is equivalent to $\lambda_n>\mu_n$
or $\lambda_n=\mu_n$ and $\lambda_{n-1}>\mu_{n-1}$ or $\cdots$.
This ordering is just opposite to $L'_n$ in page 6 of \cite{Mac}.

Let $m_0$ be a non-admissible monomial of weight 4.
Suppose that $m_0$ is of the form $e_{i_1}e_{i_2}e_{i_3}e_{i_4}$
where $0\leq i_2-i_1\leq 1, 0\leq i_4-i_3\leq 1$,
and not of the form $e_l^2e_{l+1}^2$. The relation
$\bar r_{i_1+i_2,i_3+i_4}=0$ is written as
\begin{equation}
\bar m_0=\sum_{m_k \prec m_0}c_k\bar m_k\qquad(c_k\in\mathbb{Q}).
\end{equation}
Similarly, by using $\bar r_{2l+2,2l}=0$ and $\bar r_{2l+1,2l+1}=0$,
we have
\begin{eqnarray}
&&\bar e_l^3\bar e_{l+2}=\sum_{m_k\prec\bar e_l^3\bar e_{l+2}}
c_k\bar m_k\qquad(c_k\in\mathbb{Q}),\\
&&\bar e_l^2\bar e_{l+1}^2=\sum_{m_k\prec\bar e_l^3\bar e_{l+2}}
c'_k\bar m_k\qquad(c'_k\in\mathbb{Q}).
\end{eqnarray}
Let $e_\lambda\in\pi_{n,d}$ be a non-admissible monomial.
Then there exists a non-admissible monomial $m_0$ of weight 4 such that
$e_\lambda=m_0m'$. By using the relations $\bar r_{i,j}=0$,
$\bar m_0$ can be rewritten as a linear combination of $\bar m_k$
where $m_k\prec m_0$. Hence $\bar{e}_\lambda$ is written as follows:
\begin{eqnarray}
\bar{e}_\lambda &=& \sum_{\lambda'\prec \lambda}
c_{\lambda'}\bar{e}_{\lambda'} \quad (c_{\lambda'}\in \mathbb{Q})
\end{eqnarray}
If $e_{\lambda'}$ is still non-admissible for some $\lambda'$,
we can further rewrite $\bar e_{\lambda'}$ in $R/{\cal J}$.
Since $\pi_{n,d}$ is a finite set, this procedure stops in finite times.
\end{proof}

From this proposition, we obtain an upper estimate for the dimension of
$I_{n,d}$:
\begin{equation}
\label{ineq}
\dim_{\mathbb{C}}I_{n,d}\leq\sharp\{\lambda\in
\pi_{n,d};\mbox{$\lambda$ is admissible}\}.
\end{equation}

From this follows
\begin{equation}\label{UPPER}
\dim_{\mathbb{C}(t)}J_{n,d}\leq\sharp\{\lambda\in
\pi_{n,d};\mbox{$\lambda$ is admissible}\}.
\end{equation}

\section{Construction of symmetric polynomials in $J_{n,d}$}\label{CONST}
For a partition $\lambda\in\pi_n$ and $x=(i,j)\in\lambda$, we define
\begin{eqnarray}
&&a(x)=\lambda_i-j,\ l(x)=\lambda'_j-i,\\
&&a'(x)=j-1,\ l'(x)=i-1,\\
&&n(\lambda)=\sum(i-1)\lambda_i.
\end{eqnarray}
\begin{prop}$(\hbox{\rm \cite{Mac}}$,\ {\rm Chapter VI},
$(6.11')\ ${\rm and }$(6.6))$
The specialization $u_0(P_\lambda)$ is given by
\begin{equation}\label{SPEC}
u_0(P_\lambda)=t^{n(\lambda)}\prod_{x\in\lambda}
\frac{1-t^{n-l'(x)}q^{a'(x)}}{1-t^{l(x)+1}q^{a(x)}}.
\end{equation}
We have the symmetry relation
\begin{equation}\label{sym}
\frac{u_\mu(P_\lambda)}{u_0(P_\lambda)}=\frac{u_\lambda(P_\mu)}{u_0(P_\mu)}.
\end{equation}
\end{prop}

We use the dominance ordering $\lambda>\mu$ for partitions $\lambda,\mu$.
We have
\begin{lem}\label{lem:no pole}
If $P_\lambda$ has a pole at $t^2q=1$,
 then there exists $\sigma<\lambda$ such that $E(\sigma)=E(\lambda)$ at
$t^2q=1$.
\end{lem}

From \cite{FJMM02}, we have the following proposition.

\begin{prop}\label{prop:FJMM}
$(\hbox{\rm\cite{FJMM02}}, (3.1))$ \
For $\lambda\in\pi_n$ satisfying $\lambda_i-\lambda_{i+1}\geq 2$
($1\leq i\leq n-1$), the Macdonald polynomial $P_\lambda$ has
no pole at $t^2q=1$, and the specialization $P_\lambda|_{q=t^{-2}}$
is divisible by $\prod_{1\leq i < j\leq n}(x_i-tx_j)(tx_i-x_j)$.
\end{prop}

Let $f,g\in\mathbb{C}[x_1,\ldots,x_n]\otimes_\mathbb{C}\mathbb{C}(q,t)$.
In order to prove that $f$ has no pole at $t^2q=1$,
it is sufficient to show $(1-t^2q)f=0$ at $t^2q=1$.
We take an integer $N$ such that the degree of $g$
in each variable $x_i$ is less than $N$.
In order to prove that $g=0$ at $t^2q=1$, it is sufficient to show that
there exist $n$ subsets $C_1,\ldots,C_n\subset\mathbb{C}(q,t)$,
where $\sharp(C_i)=N$, which satisfy the following two conditions:

\medskip
For each $i$ the specialization $C_i|_{q=t^{-2}}$ consists of distinct
$N$ points in $\mathbb{C}(t)$;

For all choices of $c_i\in C_i$ $(1\leq i\leq n)$,
we have $g(c_1,\cdots,c_n)=0$ at $t^2q=1$.
\begin{defn}\normalfont
Let $N$ be an integer.
A partition $\eta\in\pi_n$ is called thick if $\eta_i\gg\eta_{i+1}$ for each
$1\leq i\leq n-1$. If $\eta$ is thick, a set of $N^n$ partitions is defined
by
\[
\pi_{\eta,N}=\{\mu\in\pi_n;\mu_i=\eta_i+d_i\mbox{ for }1\leq i\leq n
\mbox{ where }0\leq d_i\leq N-1\}.
\]
For a thick partition $\eta\in\pi_{n-2}$, we define
\begin{eqnarray*}
&&\pi'_{\eta,N} = \{
\mu\in\pi_n;\mu_1=\mu_2=\eta_1+d_1,\mu_3=\mu_4=\eta_2+d_2,
\mu_i=\eta_{i-2}+d_{i-2}\\
&&\qquad\mbox{ for }5\leq i\leq n\mbox{ where }0\leq d_i\leq N-1\}.
\end{eqnarray*}
\end{defn}

We choose an sufficiently large integer $N$ and any thick partition $\eta$
when we use these sets of partitions. We do not bother to specify
$N$ and $\eta$.

\begin{defn}\normalfont
For $a\in\mathbb{C}(q,t)$, we denote by $\zeta(a)\in\mathbb{Z}$ the
multiplicity of factor $1-t^2q$ in $a$. Namely, we have
\begin{equation}
a=(1-t^2q)^{\zeta(a)}a',
\end{equation}
where the factor $a'$ has neither pole nor zero at $t^2q=1$.
\end{defn}

If $a=u_0(P_\lambda)$, using (\ref{SPEC}), we obtain
\begin{equation}\label{ZETA}
\zeta(a)=\sharp\{x\in\lambda;(l'(x),a'(x))=(n-2l,l)\}
-\sharp\{x\in\lambda;(l(x),a(x))=(2l-1,l)\}.
\end{equation}

\begin{lem}\label{lem:thickwheel}

$(i)$ Let $\eta\in\pi_n$ be a thick partition. The Macdonald polynomial
$P_\mu$
has no pole at $t^2q=1$ if $\mu\in\pi_{\eta,N}$.
Moreover, we have $\zeta(u_0(P_\mu)) = [\frac{n}{2}]$.

$(ii)$ Let $\eta\in\pi_{n-2}$ be a thick partition. The Macdonald polynomial
$P_\mu$ has no pole at $t^2q=1$ if $\mu\in\pi'_{\eta,N}$.
Moreover, we have $\zeta(u_0(P_\mu))=[\frac{n}{2}]-2$.
\end{lem}
\begin{proof}
Suppose that there exists $\sigma\in\pi_n$
such that $E(\sigma)=E(\mu)$ at $t^2q=1$.
At $t^2q=1$, $E(\mu)=\{ t^{n-i-2\mu_i} \}_{1\leq i\leq n}$. Therefore, the
condition $E(\mu)=E(\sigma)$ is equivalent to the equality of the
exponents
\[
\{n-i-2\mu_i\}_{1\leq i\leq n}=\{n-i-2\sigma_i\}_{1\leq i\leq n}
\in\mathbb{Z}^n/\mathfrak{S}_n.
\]
In the case (i), the integers $\mu_i$ are well-separated so that
we have $\mu_i=\sigma_i$ for all $i$. In the case (ii), we
can separate the indices by parity, and obtain the equalities
$\{n-2i-2\mu_{2i}\}_i=\{n-2j-2\sigma_{2j}\}_j$ and
$\{n-2i+1-2\mu_{2i-1}\}_i=\{n-2j+1-2\sigma_{2j-1}\}_j$.
{}From this we see that $\mu_i=\sigma_i$ for all $i$.
By using Lemma \ref{lem:no pole}, we conclude $P_\mu$ has no pole
in both cases.

We use (\ref{ZETA}) for the calculation of $\zeta(u_0(P_\mu))$.
In the case (i), there does not exist $x\in\mu$ such that
$(l(x),a(x))=(2l-1,l)$. On the other hand, we have
$(n-2l+1,l+1)\in\mu$ for $1\leq l\leq [\frac{n}{2}]$.
Therefore, we have $\zeta(u_0(P_\mu))=[\frac{n}{2}]$.
In the case (ii), for $x=(1,\mu_1-1),(3,\mu_3-1)\in\mu$, we have
$(l(x),a(x))=(1,1)$. Except for these two $x$'s, there does not exist
$x\in\mu$ such that $(l(x),a(x))=(2l-1,l)$. On the other hand,
we have $(n-2l+1,l+1)\in\mu$ for $1\leq l\leq[\frac{n}{2}]$.
Therefore, we have $\zeta(u_0(P_\mu))=[\frac{n}{2}]-2$.
\end{proof}

Now we are ready to prove Theorem \ref{thm:main}.
\begin{proof}[Proof of Theorem \ref{thm:main}]
Since we have (\ref{UPPER}), it is enough to show that the modified
polynomials $\bar{P}_\lambda$ satisfy the condition (\ref{WC}).

{\it Case A}.
Though we have already shown that $\bar{P}_\lambda$ has no pole
 at $t^2q=1$ in Remark \ref{rem:construction},
 we give an another proof here, which is applicable to Case B and Case C.

For $x=(i+1,\lambda_{i+1}-1)\in \lambda$ we have $(l(x),a(x))=(1,1)$,
and if $\lambda_i>2$, for $x=(i,\lambda_i-2)\in \lambda$ we have $(l(x),a(x))=(3,2)$.
Except for these two $x$'s, there does not exist $x\in\lambda$
such that $(l(x),a(x))=(2l-1,l)$. On the other hand,
we have $(n-2l+1,l+1)\in\lambda$ for $l=1,3\leq l \leq [\frac{n}{2}]$,
and if $\lambda_i>2$, $(n-3,3)\in\lambda$.
Therefore, we have $\zeta(u_0(P_\lambda))=[\frac{n}{2}]-2$.

Let $\mu\in\pi_{\eta,N}$. where $N$ is an sufficiently large integer,
and $\eta\in\pi_n$ is a thick partition. By Lemma \ref{lem:thickwheel},
$P_\mu$ has no pole at $t^2q=1$ and $\zeta(u_0(P_\mu))=[\frac{n}{2}]$.
{}From the symmetry relation (\ref{sym}), we have
\begin{eqnarray}
u_\mu(\bar{P}_\lambda)&=&u_0(P_\lambda)
\left(\frac{u_\mu(P_\lambda)}{u_0(P_\lambda)}-\frac{u_\mu(P_\nu)}{u_0(P_\nu)
}
\right)\nonumber\\
&=&u_0(P_\lambda)\frac{u_\lambda(P_\mu)-u_\nu(P_\mu)}{u_0(P_\mu)}.
\label{EQUAL}
\end{eqnarray}
Since $\mu_i-\mu_{i+1}\geq 2$, by Proposition \ref{prop:FJMM},
$P_\mu|_{q=t^{-2}}$ is divisible by
$\prod_{1\leq i < j\leq n}(x_i-tx_j)(tx_i-x_j)$.
Therefore there exists
$f_\mu,g_\mu\in\Lambda_n\otimes_\mathbb{C}\mathbb{C}(q,t)$,
 $P_\mu(x)$ is written as
\begin{eqnarray}
P_\mu(x) &=& P_\mu^{(1)}(x) + P_\mu^{(2)}(x),\\
\mbox{where, } P_\mu^{(1)}(x) &=& \prod_{1\leq i < j\leq
n}(x_i-tx_j)(tx_i-x_j)f_\mu(x)\\
\mbox{and } P_\mu^{(2)}(x) &=& (t^2q-1)g_\mu(x).
\end{eqnarray}
Since $\lambda_i=\lambda_{i+1}=\lambda_{i+2}$, we have
$u_\lambda(x_i-tx_{i+1})|_{t^2q=1}=0$ and
$u_\lambda(x_{i+1}-tx_{i+2})|_{t^2q=1}=0$.
Hence $\zeta(u_\lambda(P_\mu^{(1)}))\geq 2$.
Similarly, we have $\zeta(u_\nu(P_\mu^{(1)}))\geq 2$.
On the other hand, because $\{u_\lambda(x_i)\}_{1\leq i\leq
n}=\{u_\nu(x_i)\}_{1\leq i\leq n}$
at $t^2q=1$, we see $\zeta( u_\lambda(P_\mu^{(2)})-u_\nu(P_\mu^{(2)})
)\geq2$.
Therefore, the equality (\ref{EQUAL}) implies
$\zeta(u_\mu(\bar{P}_\lambda))\geq 0$ for all $\mu\in\pi_{\eta,N}$.
This implies that $\bar{P}_\lambda$ has no pole at $t^2q=1$.

Now, take a thick partition $\eta\in\pi_{n-2}$ and consider the set
$\pi'_{\eta,N}$. Let $\mu\in\pi'_{\eta,N}$.
By Lemma \ref{lem:thickwheel}, $P_\mu$ has no pole at $t^2q=1$
and $\zeta(u_0(P_\mu))=[\frac{n}{2}]-2$.
Since $\{u_\lambda(x_i)\}_{1\leq i\leq n}=\{u_\nu(x_i)\}_{1\leq i\leq n}$
at $t^2q=1$, we have $\zeta( u_\lambda(P_\mu)-u_\nu(P_\mu) )\geq1$.
Therefore, the equality (\ref{EQUAL}) implies
$\zeta(u_\mu(\bar{P}_\lambda))\geq 1$ for all $\mu\in\pi'_{\eta,N}$.
This implies that $\bar{P}_\lambda=0$ at $t^2q=1$.

{\it Case B}.
For $x=(i,\lambda_i-1)\in \lambda$, we have $(l(x),a(x))=(1,1)$,
and if $\lambda_{i+1}>1$, for $x=(i+1,\lambda_{i+1}-1)\in \lambda$,
 we have $(l(x),a(x))=(1,1)$.
Except for these two $x$'s, there does not exist $x\in\lambda$
such that $(l(x),a(x))=(2l-1,l)$.
On the other hand, $(n-2l+1,l+1)\in\lambda$  for $2\leq l\leq[\frac{n}{2}]$,
and if $\lambda_{i+1}>1$, $(n-1,2)\in\lambda$.
Therefore, we have $\zeta(u_0(P_\lambda))=[\frac{n}{2}]-2$.
The rest of the proof is similar to {\it Case A}.

{\it Case C}.
Take a thick partition $\eta\in\pi_n$.
Let $\mu\in\pi_{\eta,N}$. By Lemma \ref{lem:thickwheel},
$P_\mu$ has no pole at $t^2q=1$, and $\zeta(u_0(P_\mu))=[\frac{n}{2}]$.
Using (\ref{sym}), we have
\begin{equation}
u_\mu(P_\lambda)=\frac{u_\lambda(P_\mu)}{u_0(P_\mu)}u_0(P_\lambda).
\end{equation}
 There are two cases:

\noindent
(i) For all $i$, $\lambda_{i}-\lambda_{i+1}\geq 2$;

\noindent
(ii) Otherwise


For (i), we have $\zeta(u_0(P_\lambda))=[\frac{n}{2}]$.
Therefore, $\zeta(u_\mu(P_\lambda))\geq0$, and
$P_\lambda$ has no pole at $t^2q=1$.

For (ii), fix the maximum number $i_0$ such that
$\lambda_{i_0}-\lambda_{i_0+1}\leq 1$.
If $i_0=n-1$ and $\lambda_{n-1}\leq1$,
then $(i_0,\lambda_{i_0}-1)\notin\lambda$ and $(n-1,2)\notin \lambda $.
Otherwise, $(i_0,\lambda_{i_0}-1)\in \lambda$ and $(n-1,2)\in \lambda $.
In both cases, we have $\zeta(u_0(P_\lambda))=[\frac{n}{2}]-1$.
At $t^2q=1$, we have $u_\lambda(x_{i_0}-tx_{i_0+1})=0$
if $\lambda_{i_0}=\lambda_{i_0+1}$. Therefore, from Proposition
\ref{prop:FJMM}, we obtain $\zeta(u_\lambda(P_\mu))\geq 1$.
Therefore, $\zeta(u_\mu(P_\lambda))\geq 0$, and
$P_\lambda$ has no pole at $t^2q=1$.

Now, let $\eta\in\pi_{n-2}$ be a thick partition,
and $\mu\in\pi'_{\eta,N}$.
By Lemma \ref{lem:thickwheel}, $P_\mu$ has no pole at $t^2q=1$,
and $\zeta(u_0(P_\mu))=[\frac{n}{2}]-2$.
We have already shown that $\zeta(u_0(P_\lambda))\geq[\frac{n}{2}]-1$.
Therefore, $\zeta(u_\mu(P_\lambda))\geq 1$, and
$P_\lambda|_{q=t^{-2}}$ satisfies the zero condition (\ref{WC}).
\end{proof}

\section{Ideal $I_n$ and Jack polynomials with $\theta=-\frac{1}{2}$}

In this section we derive as a corollary of the previous results a linear
basis in the ideal
$I_n$ consisting of symmetric polynomials
satisfying the zero condition
\begin{equation}
f(x,x,y,y,x_5,\cdots,x_n)=0,\label{J}
\end{equation}
in terms of the Jack polynomials. Recall that the Jack polynomials are the
limiting cases of the Macdonald polynomials:
\begin{equation}
P_{\lambda}(x_1,\ldots,x_n, \theta)=\lim_{q\to1}P_{\lambda}(x_1,\ldots,x_n,
q,q^{\theta}).\label{lim}
\end{equation}
(see \cite{Mac}). We are using here the parameter $\theta =
\alpha^{-1}$ which is inverse to the parameter $\alpha$ from the
Macdonald's book \cite{Mac}. Note that the Jack polynomials are
well defined for all partitions $\lambda$ only for generic
$\theta$ (more precisely if $\theta$ is not a negative rational or
zero), so the specialization $\theta = -\frac{1}{2}$ needs some
caution.

\begin{lem}\label{lem:jack}
Consider a rational function of the form
$$
u(q,t)=\frac{f(q,t)}{\prod(1-q^kt^l)},
$$
where $f$ is a polynomial.  Suppose
that this function is well defined when $qt^2=1$ and there exists a limit
$$
\lim_{q\to 1}u(q,q^{\theta})=v(\theta)
$$
for generic $\theta.$
Then the function $v(\theta)$ is rational, well-defined at $\theta =
-\frac{1}{2}$
and
$$
\lim_{q\to 1}u(q,q^{-1/2})=v(-1/2)
$$
\end{lem}
\begin{proof}
According to our assumptions $u(q,t)$ is well defined at $qt^2=1$, so
we may suppose that $k-1/2l\ne 0$. Let us expand $f(q,q^{\theta})$
into power series using the expansion
$$
q^{\theta}=1+\theta r+\frac{\theta(\theta-1)}{2}r^2+\dots
$$
where $r=q-1:$
$$
f(q,q^{\theta})=\sum_{i=0}^{\infty}a_{i}(\theta)r^i
$$
where $a_i(\theta)$ are polynomials. Let $p$ be the number of
factors in the denominator of $u$.  From our assumptions we have
$$
a_{0}=a_{1}=\dots=a_{p-1}=0
$$
and
$$
\lim_{q\to1}u(q,q^{\theta})=\frac{a_{p}}{\prod(k+\theta l)},
$$
which implies the lemma.
\end{proof}

Now we can define the {\it modified Jack polynomials} $\bar
P_{\lambda}(x,\theta)$ for the admissible partitions of type A and
B (for type C it is just usual Jack polynomial) by the formula
(\ref{BAR}) with $P_\lambda$ replaced by $P_{\lambda} (x,\theta)$ and
$u_0(P_\lambda)$ by
\begin{equation}
\lim_{q\to1}u_{0}(P_{\lambda}(x,q,q^{\theta}))=P_{\lambda}(1,\dots,1,\theta)
=
\prod_{s\in\lambda}\frac{(n-l^{\prime}(s))\theta+a^{\prime}(s)}
{(l(s)+1)\theta+a(s)}.\label{J1}
\end{equation}
The last relation follows from the formula (\ref{SPEC}).
We denote this quantity by $u_0(P_{\lambda}(x,\theta))$.

\begin{thm}\label{thm:main1}
For any admissible diagram $\lambda$ the modified Jack polynomial
$\bar{P}_{\lambda}(x,\theta)$ is well defined at $\theta=-1/2$ and
the specialization $\bar{P}_{\lambda}(x,\theta)|_{\theta=-1/2}$
satisfies the zero condition $(\ref{J})$. The set of polynomials
$\{\bar{P}_{\lambda}(x,\theta)|_{\theta=-1/2};\lambda\in\pi_{n}
\mbox{ is admissible}\}$ is a linear basis of the ideal $I_n$.
\end{thm}
\begin{proof}
From the Theorem \ref{thm:main} and Lemma \ref{lem:jack} we see that
$\bar{P}_{\lambda}(x,\theta)$ are well defined at $\theta =-1/2$
and satisfy the condition (\ref{J}). It is easy to check that they are
linearly independent.
Now the theorem follows from the inequality (\ref{ineq}) .
\end{proof}

\section{ Ideal $I^*_n$ and interpolation Jack polynomials}

Consider now the ideal  $I^*_{n}$ of symmetric polynomials which
satisfy the following zero condition
\begin{equation}
f(x, x+1/2, y, y+1/2, x_5,\cdots, x_n)=0.\label{IJ}
\end{equation}
We will construct a basis of $I^*_{n}$ by using interpolation Jack
polynomials specialized at $\theta=-1/2$.

Let us assume first that $\theta$ is not negative rational
or zero and denote by $\rho(\theta)$ the following ordered set
$((n-1)\theta,(n-2)\theta,\dots,\theta,0)$.

\begin{defn}\normalfont (\cite{KS,OO})
Let $\lambda$ be a partition with $\lambda_{n+1}=0$. There exists a
unique symmetric polynomial $P^*(x,\theta)$ such that

1)\quad $deg P_{\lambda}^*(x,\theta)\le |\lambda|$

2)\quad $P_{\lambda}^*(\mu+\rho(\theta),\theta)=0$ \quad if \quad
$|\mu|\le |\lambda|$, $\mu\ne\lambda,\mu_{n+1}=0$

3)\quad $P_{\lambda}^*(\lambda+\rho(\theta),\theta)=
\prod_{s\in\lambda}(l(s)\theta+a(s)+1)$

\end{defn}
These polynomials are called {\it interpolation Jack polynomials} (cf.
\cite{O}).

Let us define now for admissible partitions of type A and B the
{\it modified interpolation Jack polynomial} $\bar
P_{\lambda}^*(x,\theta)$ by the formula (\ref{BAR}), where by
definition
$u_{0}(P_{\lambda}^*(x,\theta))=u_{0}(P_{\lambda}(x,\theta)).$ For
partitions of type C we have just usual interpolation Jack
polynomials.

\begin{thm}\label{thm:main2}
For any admissible partition the modified polynomial
$\bar{P}^*_{\lambda}(x,\theta)$ is well-defined at \,
$\theta=-1/2$ and the specialization
$\bar{P}_{\lambda}(x,\theta)^*|_{\theta=-1/2}$ satisfies the zero
condition $(\ref{IJ})$. The set of polynomials
$\{\bar{P}^*_{\lambda}(x,\theta)|_{\theta=-1/2};\lambda\in\pi_{n}
\mbox{ is admissible}\}$ is a basis of $I^*_n$.
\end{thm}
\begin{proof}

The zero condition (\ref {IJ}) is equivalent to the following
relation
\begin{equation}
f(\mu+ \rho(-\frac{1}{2})) = f( \mu_{1}-\frac{n-1}{2},\,
\mu_{2}-\frac{n-2}{2},\dots,\,\mu_{n-1}-\frac{1}{2}, \mu_n ) =
0\label{IJ2}
\end{equation}
for any non-admissible partition $\mu$.

To prove the relation (\ref{IJ2}) for the modified interpolation
Jack polynomials we need the following variant of the Pieri
formula:
\begin{equation}
\left(\sum_{i=1}^n
x_{i}-|\lambda|\right)\bar{P}^*_{\lambda}(x,\theta)|_{\theta=-1/2}=
\sum_{\tau}c_{\lambda,\tau}\bar{P}^*_{\tau}(x,\theta)|_{\theta=-1/2},
\label{IJ3}
\end{equation}
 where $\lambda$ is an admissible partition and the sum is taken over
admissible $\tau$ such that $|\tau|=|\lambda|+1$.

To prove (\ref{IJ3}) note that for generic $\theta$ the both sets
$P_{\lambda}(x,\theta)$ and $P_{\lambda}^*(x,\theta)$ are the
bases in the algebra $\Lambda_n$ of symmetric polynomials in
$x_{1},\dots,x_{n}$, so there exists a linear isomorphism
$\Psi_{\theta}: \Lambda_n \rightarrow \Lambda_n$ such that $$
\Psi_{\theta}(P_{\lambda}(x,\theta))=P_{\lambda}^*(x,\theta) $$ We
will call it the {\it dehomogeneization operator.} Knop and Sahi
in \cite{KS} found the following explicit formula for this
operator.

Let us define difference operators
${\cal{E}}_{1},\dots,{\cal{E}}_{n}$ from the following equality
$$
1+{\cal{E}}_{1}t+\dots+{\cal{E}}_{n}t^n=\prod_{i<j}(x_{i}-x_{j})^{-1}
det\left[(x_{i}+\theta)^{n-j}+tx_{i}^{n-j+1}T_{i}\right],
$$
where $T_{i}$ is the shift operator:
$$
(T_{i}f)(x_{1},\dots,x_{i},\dots,x_{n})=f(x_{1},\dots,x_{i}-1,\dots,x_{n}).
$$
One can show that the operators ${\cal{E}}_{1},\dots,{\cal{E}}_{n}$ commute
with each
other.

Let $\psi_{\theta}: \Lambda_{n}\longrightarrow {\bf C}
[{\cal{E}}_{1},\dots,{\cal{E}}_{n}]$ be a linear isomorphism,
sending the $k$-th elementary symmetric function into
${\cal{E}}_{k},$ then
\begin{equation}
\Psi_{\theta}(f)=\psi_{\theta}(f)(1). \label{IJ1}
\end{equation}
(see \cite{KS}).
Now we are ready to prove
(\ref{IJ3}). Since $I_{n}$ is an ideal we have
\begin{equation}
\left(\sum_{i=1}^n
x_{i}\right)\bar{P}_{\lambda}(x,\theta)|_{\theta=-1/2}=
\sum_{\tau}c_{\lambda,\tau}\bar{P}_{\tau}(x,\theta)|_{\theta=-1/2},\label{IJ4}
\end{equation}
where $\lambda$ is an admissible partition and the sum is taken
over the admissible $\tau$ such that $|\tau|=|\lambda|+1$. Now if
we apply to both sides of (\ref{IJ4}) the dehomogeneization
operator $\Psi_{-1/2}$ and take into account the relation
\cite{KS} $$ \Psi_{\theta}\left(\sum_{i=1}^n x_{i}f\right
)=\left(\sum_{i=1}^n x_{i}-\deg f \right)\Psi_{\theta}(f) $$
(assuming that $f$ is homogeneous) we get (\ref{IJ3}).

To prove the vanishing property (\ref{IJ2}) for the modified
interpolation Jack polynomials we first show that
\begin{equation}
 \bar{P}^*_{\lambda}(\mu+\rho(-1/2),\theta)|_{\theta=-1/2}=0\label{IJ5}
\end{equation}
for any non-admissible $\mu$ such that
$|\mu|\le|\lambda|$ .

Consider first the case A. If $\mu\ne\nu$ this
follows from the definition.  If $\mu=\nu$
then we have
$$
\bar{P}^*_{\lambda}(\nu+\rho(\theta),\theta)={P}^*_{\lambda}(\nu+\rho(\theta
),\theta)-
\frac{u_{0}(P^*_{\lambda}(x,\theta))}
{u_{0}(P^*_{\nu}(x,\theta))}{P}^*_{\nu}(\nu+\rho(\theta),\theta)=
$$
$$
-\frac{u_{0}(P^*_{\lambda}(x,\theta))}
{u_{0}(P^*_{\nu}(x,\theta))}{P}^*_{\nu}(\nu+\rho(\theta),\theta).
$$
For $f \in\mathbb{C}(\theta)$ we denote by $\zeta(f)\in \mathbb{Z}$
the multiplicity of factor $\theta+1/2$ in $f$. It is not
difficult to verify that
\begin{equation}
\label{last1}
\zeta(u_{0}(P^*_{\lambda}(x,\theta)))=
\zeta(u_{0}(P^*_{\nu}(x,\theta)))=[\frac{n}{2}]-2
\end{equation}
and
\begin{equation}
\label{last2}
\zeta({P}^*_{\nu}(\nu+\rho(\theta),\theta))>0.
\end{equation}

The formula (\ref{last1}) can be proved in the same way as before (see
section 4).
To prove (\ref{last2}) let us note that from the property 3 in the
Definition 6.1 it 
follows that ${P}^*_{\nu}(\nu+\rho(\theta),\theta)$ is a polynomial in
$\theta$ which has a zero at $\theta = -1/2$ whenever the diagram of $\nu$
has a box $s$ with $(l(s),a(s))=(2l,l-1)$.
Since for a partition $\lambda$ with $\lambda_i=\lambda_{i+1}=\lambda_{i+2}
= \lambda_{i+3}+2$ the corresponding diagram $\nu$ has a box
$s=(i+1,\nu_{i+1})$ with
$l(s)=2, a(s)=0$ this implies (\ref{last2}). This  proves (\ref{IJ5}) in the
case $A$. 

In the cases B and C the relation (\ref{IJ5})
follows from the definition of the interpolation Jack polynomials.

Now we use the formula (\ref{IJ3}) and the induction in $|\mu|$ to prove
that the relation (\ref{IJ5}) is valid for all non-admissible $\mu$. The
theorem now follows.
\end{proof}

\section{Some corollaries and discussion.}

In this section, first, we discuss the character
(or the Hilbert-Poincar\'e series) for the ideal $I_n$
\begin{eqnarray}
{\mathop{\mathrm{ch}}\nolimits}\,
I_n:=\sum_{d=0}^{\infty} (\dim I_{n,d}) q^d \nonumber
\end{eqnarray}
We denote the product $\prod_{j=1}^{i}(1-q^j)$ by $(q)_i$.
We have already constructed the basis of $I_n$ which is labeled
by admissible partitions (Theorem \ref{thm:main1}).
Let us count admissible partitions $\lambda$.
There are four different cases:

\noindent
{\it Case 1}
\quad For any $i$,
\begin{eqnarray}
&&\lambda_i-\lambda_{i+1}\geq 2 \nonumber
\end{eqnarray}
In this case, the character of such partitions is given by
\begin{eqnarray}
\frac{q^{n(n-1)}}{(q)_n}. \label{char1}
\end{eqnarray}
{\it Case 2}
\quad There exists $i$ such that
\begin{eqnarray}
&&\lambda_i=\lambda_{i+1}=\lambda_{i+2}, \nonumber\\
&&\lambda_{i-1}-\lambda_i\geq 3 \quad\mbox{if $i\geq 2$}, \nonumber\\
&&\lambda_j-\lambda_{j+1}\geq 2 \quad\mbox{if $j\leq i-2$ or $j\geq i+2$}. \nonumber
\end{eqnarray}
The character of this case is
\begin{eqnarray}
\frac{q^{(n-1)(n-3)}}{(q)_n}\sum_{i=0}^{n-3}q^{3i}(1-q^{n-2-i})(1-q^{n-1-i}). \label{char2}
\end{eqnarray}
{\it Case 3}
\quad There exists $i$ such that
\begin{eqnarray}
&&\lambda_i=\lambda_{i+1}, \nonumber\\
&&\lambda_{i-1}-\lambda_i\geq 1 \quad\mbox{if $i\geq 2$}, \nonumber\\
&&\lambda_j-\lambda_{j+1}\geq 2 \quad\mbox{if $j\leq i-2$ or $j\geq i+1$}. \nonumber
\end{eqnarray}
The character of this case is
\begin{eqnarray}
\frac{q^{(n-2)^2}}{(q)_n}\sum_{i=0}^{n-2}q^{3i}(1-q^{n-1-i}). \label{char3}
\end{eqnarray}
{\it Case 4}
\quad There exists $i$ such that
\begin{eqnarray}
&&\lambda_i=\lambda_{i+1}+1, \nonumber\\
&&\lambda_{i-1}-\lambda_i\geq 0 \quad\mbox{if $i\geq 2$}, \nonumber\\
&&\lambda_j-\lambda_{j+1}\geq 2 \quad\mbox{if $j\leq i-2$ or $j\geq i+1$}. \nonumber
\end{eqnarray}
The character of this case is
\begin{eqnarray}
\frac{q^{(n-2)^2+1}}{(q)_n}\sum_{i=0}^{n-2}q^{3i}(1-q^{n-1-i}). \label{char4}
\end{eqnarray}
\medskip
\noindent
Summing up these characters, we obtain
\begin{thm}\label{PCH}
The character for the ideal $I_n$ is given by
\begin{eqnarray}
{\mathop{\mathrm{ch}}\nolimits}\,I_n=\frac{q^{(n-1)(n-3)}}{(q)_1(q)_{n-3}}
+\frac{q^{(n-2)(n-1)}}{(q)_1(q)_{n-2}}
+\frac{q^{n(n-1)}}{(q)_n}. \label{char}
\end{eqnarray}
\end{thm}

We give an interpretation of the expression (\ref{char}).
Define
\begin{eqnarray}
F_2&:=&\{f\in\Lambda_n ; f(x,x,y,y,x_5,\cdots,x_n)=0 \}, \nonumber\\
F&:=&\left\{f\in\Lambda_n ; \begin{array}{l} f(x,x,y,y,x_5,\cdots,x_n)=0 \mbox{ and } \\ f(x,x,x,x_4,\cdots,x_n)=0 \end{array} \right\}, \nonumber\\
F_1&:=&\{f\in\Lambda_n ; f(x,x,x_3,\cdots,x_n)=0 \}. \nonumber
\end{eqnarray}
Note that $F_2 \supset F \supset F_1$.
Then, we have injective maps:
\begin{eqnarray}
F_2/F\ni\bar{f}
&\mapsto& \prod_{i=4}^{n}(z-x_i)^3 \prod_{4\leq i<j\leq n}(x_i-x_j)^2
g(z,x_4,\cdots,x_n) \nonumber\\
&&\hbox{ where }
g\in\mathbb{C}[z]\otimes\mathbb{C}[x_4,\cdots,x_n]^{\mathfrak{S}_{n-3}},
\nonumber\\
F/F_1\ni\bar{f} &\mapsto& \prod_{i=3}^{n}(z-x_i)^2 \prod_{3\leq i<j\leq n}(x_i-x_j)^2 h(z,x_3,\cdots,x_n)\nonumber\\&&\hbox{ where }
h\in\mathbb{C}[z]\otimes\mathbb{C}[x_3,\cdots,x_n]^{\mathfrak{S}_{n-2}}),
\nonumber\\
F_1\ni f &\mapsto& \prod_{1\leq i<j\leq n}(x_i-x_j)^2 k(x_1,\cdots,x_n)
\nonumber\\
&&\hbox{ where }k\in\mathbb{C}[x_1,\cdots,x_n]^{\mathfrak{S}_{n}}). \nonumber
\end{eqnarray}
From injectivity of these maps, thier characters are estimated from above:
\begin{eqnarray}
{\mathop{\mathrm{ch}}\nolimits}\,F_2/F\leq
\frac{q^{(n-1)(n-3)}}{(q)_1(q)_{n-3}},  \quad
{\mathop{\mathrm{ch}}\nolimits}\,F/F_1\leq
\frac{q^{(n-2)(n-1)}}{(q)_1(q)_{n-2}},  \quad
{\mathop{\mathrm{ch}}\nolimits}\,F_1\leq\frac{q^{n(n-1)}}{(q)_n}. \nonumber
\end{eqnarray}
Theorem  \ref{PCH} imples that
the equalities hold in the above. Namely, the above
mappings are, in fact, isomorphisms.

From algebraic point of view a more natural question is what are
the algebraic generators of an ideal rather than what is its
linear basis. As a corollary of our results we see that the minimal degree of the generators in the ideal $I_n$ is
$$M(n)=(n-1)(n-3) = n^2 - 4n + 3.$$ The corresponding generator is
given by the Jack polynomial $P_{\lambda}(x,\theta)$ with $\theta
= -\frac{1}{2}$ and the partition $\lambda = ((2n-5), (2n-7),
\ldots, 5, 3, 0, 0, 0)$,
which has the minimal
weight $|\lambda|$ among all the admissible partitions with $n$
parts (if $n$ is larger than $3$).
It is easy to see that the polynomial $$Q(x) = Symm \prod_{j=4}^{n} (x_1 - x_j)(x_2 - x_j)(x_3 - x_j)
\prod_{k,l=4,\, k <l}^{n}(x_k - x_l)^2,$$
where $Symm$ means symmetrization, 
belongs to the ideal $I_n$ and has the degree $M(n)$, and thus must coincide up to a multiple with this Jack polynomial.
 For $n=4$ and $n=5$ we have
respectively $M=3$ and $M=8$ in agreement with Cayley's results
\cite{Cayley} and with Magma calculations kindly performed for us
by Miles Reid: the total set of the degrees of the generators are

$n=4$ \quad \quad 3, 4, 5, 6, 7, 8, 9

$n=5$ \quad \quad 8, 9, 10, 11, 12, 13, 14, 15, 16, 17

$n=6$ \quad \quad 15, 16, 17, 18, 19, 20, 21, 22, 23, 24, 25, 26,
27.

\noindent The computer calculations for larger $n$ need
substantially more time to perform, so it is not easy to check if
this very suggestive pattern holds for higher degrees. But at
least we see that the Cohen-Macaulay property in general is not
satisfied.

One can consider a general question about coincident root loci
related to any partition $\mu$ of the degree of the polynomial
$n$ when we have the roots of the multiplicities $\mu_1,
\mu_2, \dots, \mu_k.$ Our problem corresponds to the
partitions $(2,2,1,1,\dots,1)= (2^2,1^{(n-4)}).$ The results of
the paper \cite{SV} suggest that we should expect similar results
for the partitions $(p^k, q^l),\, k p + l q = n.$ A very
interesting question is what happens for other partitions.

For the ideals $I_n(p)$ corresponding to the partitions $(p,
1^{n-p})$ (describing the polynomials of degree $n$ with a root of
multiplicity at least $p$) a formula for the Hilbert-Poincar\`e series was given 
by Feigin and Stoyanovsky in \cite{FS} and 
a linear basis in terms of Jack
polynomials was constructed in the paper \cite{FJMM}. 
A simple analysis of the admissible diagrams from \cite{FJMM} leads to the
following formula for the minimal degree of the generators in the ideal $I_n(p):$ 
$$M(n,p)= s(s-1)(p-1)+ 2sr,$$
where $s$ is the result of the division of $n$ by $p-1$ and $r$ is
the corresponding remainder: $n = s(p-1) + r.$
The corresponding generator is given
by the Jack polynomial $P_{\lambda}(x,\theta)$ with $\theta =
-\frac{1}{p-1}$ and the partition $\lambda =
((2s)^{r},(2s-2)^{p-1}, (2s-4)^{p-1}, \ldots, 2^{p-1}, 0^{p-1}).$

The first non-trivial case $p=3$ corresponds to the polynomials
with a triple root. In that case the minimal degree of the
generators depends on the parity of $n$ and equals to the integer
part of $\frac{1}{2}(n-1)^2.$ This is in a good agreement with
Magma calculation of the degrees of the generators for $n \leq 8$
(calculations for $n$ bigger than $9$ seem to take unreasonable
computer time to complete):

$n=4$ \quad \quad 4, 6

$n=5$ \quad \quad 8, 9, 10, 10, 12

$n=6$ \quad \quad 12, 14, 15, 16, 18

$n=7$ \quad \quad 18, 19, 20, 20, 21, 21, 22, 22, 23, 24

$n=8$ \quad \quad 24, 26, 27, 28, 28, 29, 30, 30, 31, 32, 32, 33,
34

The best results about the algebraic generators for the partitions
$(p, 1^{n-p})$ we are aware of belong to J. Weyman \cite{W1,W2}.
In particular he proved that if $p$ is larger than the half of $n$
the corresponding ideal is generated by the polynomials of degree
4 (in the coefficients) and made a conjecture about the degrees of
the generators in general case. We should note that since we are
using different grading the comparison with these results is not
straightforward.

It is clear that the relations between Sylvester-Cayley problem
and the theory of Jack and Macdonald polynomials with special
values of the parameters should be investigated further. In
particular it would be interesting to understand how the classical
invariant theory is related to the representation theory behind
these polynomials.

{\it Acknowledgements.}
We are grateful to the organisers of Jack-Macdonald meeting at
ICMS (Edinburgh, September 2003) where this work was initiated. We
would like to thank Andrei Okounkov for useful comments and
especially Miles Reid for providing us with the results of
computer algebra experiments and many stimulating discussions.

\end{document}